\documentclass[11pt]{article}
\usepackage{}
\usepackage[total={6in, 9in}]{geometry}
 \usepackage{amsfonts}
\usepackage{amsthm}
\usepackage{amssymb}
\usepackage{amsmath,amsfonts}
\usepackage{mathrsfs}
\usepackage{cases}
\usepackage{latexsym,bm}
\usepackage{indentfirst}
\usepackage{color}
\usepackage{ifpdf}
\usepackage{graphicx}
\usepackage{psfrag}
\usepackage{caption}
\usepackage{enumerate}
\usepackage{dsfont}
\usepackage[
pdfauthor={},
pdftitle={},
pdfstartview=XYZ,
bookmarks=true,
colorlinks=true,
linkcolor=blue,
urlcolor=blue,
citecolor=blue,
bookmarks=true,
linktocpage=true,
hyperindex=true
]{hyperref}

\usepackage[natural]{xcolor}
\usepackage{tikz}
\usepackage{subfigure,tikz}
\usepackage{tikz-3dplot}
\usepackage{pgf,tikz,pgfplots}
\usepackage{mathrsfs}

\newtheorem{THM}{\textbf{Theorem}}

\newtheorem{CLA}{\textbf{Claim}}[section]
\newtheorem{CON}[THM]{\textbf{Conjecture}}
\newtheorem{COR}[THM]{\textbf{Corollary}}

\newcommand{\qqed}{\hfill $\blacksquare$\vspace{1mm}}

\newcommand{\iP}{\overset{\leftharpoonup }{P}}

\newcommand{\oC}{\overset{\rightharpoonup }{C}}
\newcommand{\oP}{\overset{\rightharpoonup }{P}}

\newcommand{\lC}{\overset{\leftharpoonup}{C}}
\newcommand{\rC}{\overset{\rightharpoonup }{C}}
\newcommand{\lP}{\overset{\leftharpoonup }{P}}
\newcommand{\rP}{\overset{\rightharpoonup }{P}}

\begin{document}
\title{Degree sequence condition for Hamiltonicity  in tough graphs}
\author{Songling Shan\footnote{Auburn University, Department of Mathematics and Statistics, Auburn, AL 36849.
		Email: {\tt szs0398@auburn.edu}.   
		Supported in part by  NSF grants DMS-2345869 and DMS-2451895.}
		\qquad 
		 Arthur Tanyel\footnote{Auburn University, Department of Mathematics and Statistics, Auburn, AL 36849.
		 	Email:	{\tt ant0034@auburn.edu}. }
}

\date{\today}
\maketitle

\begin{abstract}
Generalizing both Dirac's condition and Ore's condition for Hamilton cycles, Chv\'atal in 1972 established a degree sequence condition
for the existence of a Hamilton cycle in a graph.  Ho\`ang in 1995 generalized Chv\'atal's degree sequence condition for 1-tough graphs and 
conjectured  a $t$-tough analogue for any positive integer $t\ge 1$. Ho\`ang in the same paper verified his conjecture for $t\le 3$
and recently Ho\`ang and Robin verified the conjecture for $t=4$. In this paper, we confirm the conjecture for all $t\ge 4$. 
The proof depends on two newly established  results on cycle structures in tough graphs, which hold independent interest. 

\medskip 


\end{abstract}

\section{Introduction}

Graphs considered in this paper are simple, undirected, and finite.
Let $G$ be a graph.
Denote by $V(G)$ and  $E(G)$ the vertex set and edge set of $G$,
respectively.  For a vertex $v\in V(G)$, the set of neighbors of $v$ is denoted by $N(v)$, and the  degree  of $v$  is denoted 
by $\deg(v)$, which is the size of $N(v)$. If $H\subseteq G$, then $\deg(v, H)$ is the size of the set $N(v)\cap V(H)$.  If $u$ and $v$ are non-adjacent in $G$, then $G+uv$ is 
obtained from $G$ by adding the edge $uv$. We write $u\sim v$ if two vertices $u$ and $v$
are adjacent in $G$ and write $u\not\sim v$ otherwise.   For	 $S \subseteq V(G)$, denote by 
$G[S]$ and $G-S$ the subgraph of $G$ induced on $S$ and $V(G)\setminus S$, respectively. 
For $v\in V(G)$,  we write $G-v$ for $G-\{v\}$.
For two integers, $p$ and $q$, we let $[p,q]=\{i\in \mathbb{Z}:  p\le i \le q\}$. 

Let $n\ge 1$ be an integer. 
The non-decreasing sequence $d_1, d_2, \dots , d_n$ is a \textit{degree sequence} of a  graph $G$ if the vertices of $G$ can be labeled as $v_1, v_2, \dots, v_n$ such that $\deg(v_i) = d_i$ for all $i\in [1,n]$. In 1972, 
Chv\'atal~\cite{chvatal1972hamilton} proved the following well known result. 
\begin{THM}\label{thm:chv_degree}
	Let $G$ be a graph with degree sequence $d_1, d_2, \dots, d_n$, where $n\ge 3$ is an integer. If for all $i < \frac{n}{2}$, $d_i \le i$ implies $d_{n-i} \ge n-i$, then $G$ is Hamiltonian.
\end{THM}

Ho\`ang~\cite[Conjecture 1]{hoang1995hamiltonian} in 1995 conjectured a toughness analogue for the theorem above. 
We let $c(G)$ denote the number of components of $G$. For a real number $t \ge 0$, we say $G$ is \textit{t-tough} if $|S| \ge t \cdot c(G - S)$ for all $S \subseteq V(G)$ such that $c(G - S) \ge 2$. The largest $t$ for which $G$ is $t$-tough is called the \textit{toughness} of $G$ and is denoted $\tau(G)$. If $G$ is complete,  then $\tau (G)$ is defined to be $\infty$. Chv\'atal~\cite{chvatal1973tough} defined this concept in 1973 as a measure of a graph's ``resilience" under the removal of vertices.  
Ho\`ang's conjecture can now be stated as follows. 
\begin{CON}\label{con:degree-sequence-t-tough}
	Let  $n\ge 3$ and $t\ge 1$ be integers, and $G$ be  a $t$-tough  graph with degree sequence $d_1, d_2, \dots , d_n$.  If  for all $i < \frac{n}{2}$ it holds that  $d_i \le i$ implies  $d_{n-i+t} \ge n-i$, then  $G$ is Hamiltonian.
\end{CON}

Ho\`ang in the same paper~\cite[Theorem 3]{hoang1995hamiltonian} proved the conjecture  for $t \le 3$. 
Since every Hamiltonian graph must necessarily be 1-tough, the statement for $t=1$ generalizes Theorem~\ref{thm:chv_degree}. 
Recently, Ho\`ang  and Robin~\cite{hoang2024closure} proved the conjecture  for $t = 4$.
In this paper, we  confirm  Conjecture~\ref{con:degree-sequence-t-tough}  for all $t \ge 4$.  

\begin{THM}\label{main}
	Let $t \ge 4$ be an integer and  $G$ be  a $t$-tough graph on $n \geq 3$ vertices with degree sequence $d_1, d_2, \dots , d_n$.  If  for all $i < \frac{n}{2}$ it holds that  $d_i \le i$ implies  $d_{n-i+t} \ge n-i$, then  $G$ is Hamiltonian.
\end{THM}

A graph $G$ is \emph{pancyclic} if $G$ contains cycles of any length from $3$ to $|V(G)|$. 
As a consequence of Theorem~\ref{main}, 
a result of Ho\`ang~\cite[Theorem 7]{hoang1995hamiltonian} that if a $t$-tough graph $G$ satisfies the degree sequence condition in Theorem~\ref{main} and is Hamiltonian, then $G$ is pancyclic or bipartite,  and  the fact that bipartite graphs of order at least three have toughness at most one, 
 we also obtain the following result. 

\begin{COR}\label{cor:pancyclic}
Let $t \ge 4$ be an integer and  $G$ be  a $t$-tough graph on $n \geq 3$ vertices with degree sequence $d_1, d_2, \dots , d_n$.  If  for all $i < \frac{n}{2}$ it holds that  $d_i \le i$ implies  $d_{n-i+t} \ge n-i$,  then $G$ is pancyclic. 
\end{COR}

The proof of Theorem~\ref{main} relies on  our  closure lemma for $t$-tough graphs $G$ as stated below. 

\begin{THM}[Toughness Closure Lemma]\label{thm:t-closure}
	Let $t \ge 4$ be a rational number,   $G$ be  a $t$-tough graph on $n \geq 3$ vertices, and let   distinct  $x, y \in V(G)$  be non-adjacent with $\deg(x) + \deg(y) \geq n - t$. 
	Then $G$ is Hamiltonian if and only if  $G+xy$ is Hamiltonian.  
\end{THM}

The proof of Theorem~\ref{thm:t-closure}  relies on a restricted cycle structure in tough graphs as stated in Theorem~\ref{thm:sum-of-successors}. 
We define some notation in order to state the theorem.

Denote by $\oC$ an orientation of a cycle $C$.
We assume that the orientation is clockwise throughout the rest of this paper.
For $u,v \in V(C)$, $u \rC v$ denotes the path from $u$ to $v$ along $\oC$.  Similarly, $u \lC v$ denotes the path between $u$ and $v$ which travels opposite to the orientation. 
We use $u^+$ to denote the immediate successor of $u$ on $\oC$  and $u^-$ to denote the immediate predecessor of $u$ on $\oC$.
If $S \subseteq V(C)$, then $S^+ = \{u^+:  u \in S\}$ and $S^- = \{u^-: u \in S\}$. 
We use similar notation for a path $P$ when it is given an orientation.  

\begin{THM}\label{thm:sum-of-successors}
	Let $t\ge 4$ be rational and  $G$ be  a $t$-tough graph on $n \geq 3$ vertices. Suppose that $G$ is not Hamiltonian, but there exists $z \in V(G)$ such that $G -z$ has a Hamilton cycle $C$. Then, for any distinct $x, y \in (N(z))^+$, we have that $\deg(x)+\deg(y) < n-t$. 
\end{THM}

We will prove Theorem~\ref{main} in the next section by applying  Theorem~\ref{thm:t-closure}. 
Then Theorems~\ref{thm:t-closure} and~\ref{thm:sum-of-successors} are respectively proved in 
Sections 3 and 4.

\section{Proof of Theorem~\ref{main}}

We will need the following result by Bauer et al.~\cite{bauer1995long} and our closure lemma for $t$-tough graphs with $t\ge 4$. 

\begin{THM}\label{bauer}
	Let  $t \ge 0$ be any real number and   $G$ be a $t$-tough graph on $n \ge 3$ vertices. If  $\delta(G) > \frac{n}{t+1} - 1$, then  $G$ is Hamiltonian.
\end{THM}

%\begin{THM}[Toughness Closure Lemma]\label{thm:t-closure}
%	Let $t \geq 4$ be an integer,   $G$ be  a $t$-tough graph on $n \geq 3$ vertices, and let   distinct  $x, y \in V(G)$  be non-adjacent with $\deg(x) + \deg(y) \geq n - t$. 
%	Then $G$ is Hamiltonian if and only if  $G+xy$ is Hamiltonian.  
%\end{THM}

The following toughness closure concept was given by Ho\`ang and Robin~\cite{hoang2024closure}.  Let $t \geq 1$ be an integer,  and $G$ be  a $t$-tough graph on $n \geq 3$ vertices. 
Then the $t$-closure   of $G$ is formed  by repeatedly adding edges joining vertices $x$ and $y$ such that $x$ and $y$ are non-adjacent in the current graph and
their degree sum is at least $n-t$ in the current graph, until no such pair remains.  By the same argument that shows that the Hamiltonian closure  of a graph is well defined (e.g., see~\cite[Lemma 4.4.2]{Bonday-and-Murty-Book}),  the 
$t$-closure   of $G$  is well defined.  Thus by Theorem~\ref{thm:t-closure},  we will consider the $t$-closure of $G$ instead of $G$ itself when we prove Theorem~\ref{main}. 
Once the closure lemma is established, the proof of Theorem~\ref{main} follows a standard argument, akin to that in Ho\`ang and Robin's work~\cite{hoang2024closure}.

\proof[Proof of Theorem~\ref{main}] 

Because $G$ satisfying  the degree sequence condition implies that any supergraph of $G$ obtained from $G$ by adding missing edges also satisfies the  degree sequence condition, 
by Theorem~\ref{thm:t-closure}, it suffices to work with the $t$-closure of $G$. For the sake of notation, we just assume that $G$ itself is its $t$-closure.   We may assume that $G$ is not Hamiltonian. Thus $G$ is not complete 
and so $\delta(G) \ge 8$ by $G$ being 4-tough.

Let $v_1, v_2, \dots, v_n$  be all the vertices of $G$ such that $\deg(v_i) = d_i$ for all $i\in [1,n]$.  
Thus, we have that $\deg(v_i) + \deg(v_j) \ge n - t$ implies $v_iv_j \in E(G)$. By Theorem~\ref{thm:chv_degree}, if $d_i > i $ for all $i < \frac{n}{2}$, then $G$ is Hamiltonian. So, we assume  that 
there exists some positive  integer $k < \frac{n}{2}$ such that $d_k \leq k$.  Then as  $\delta(G) \ge 8$, we have $k\ge 8$. 
Choose $k$ to be minimum  with the property that $d_k \leq k$.  Then  $d_i > i$ for all $i\in[1,k-1]$.  Since $d_{k-1} \le d_k \le k$, we must have $d_{k-1} = d_k = k$. 

Let $S,T \subseteq V(G)$. We say that $S$ is \textit{complete to} $T$ if for all $u \in S$ and $v \in T$ such that $u \ne v$,  we have $u \sim v$. If $u \sim v$ for all $u \in S$ and $v \in V(G)$ such that $u \ne v$, we call $S$ a \textit{universal clique} of $G$.  Clearly, 
vertices in a universal clique have degree $n-1$ in $G$. 
We will show that $G$ has a universal clique of size larger than $\frac{n}{t+1}-1$. In particular, this gives that $\delta(G) > \frac{n}{t+1}-1$. By Theorem~\ref{bauer}, this proves that $G$ is Hamiltonian, a contradiction to the assumption that $G$ is not Hamiltonian.  Let 
$$
U^\alpha = \{v_i :  d_i \ge n - \alpha, i\in [1,n]\} \quad \text{for any  integer $\alpha$ with $1\le \alpha <\frac{n}{2}$.}
$$

\medskip

\begin{CLA}\label{4.1}  For all  positive integers $\alpha < \frac{n}{2}$, $U^\alpha $ is a clique complete to $\{ v_i:  d_i \ge \alpha - t, i\in [1,n]\}$.
\end{CLA}

\proof[Proof of Claim~\ref{4.1}] If $v_j \in U^\alpha$  for some $j\in [1,n]$ and $v_\ell \in \{ v_i :  d_i \ge \alpha - t, i\in[1,n]\}$ for some $\ell\in [1,n]$, then $d_j + d_\ell \ge n - \alpha+\alpha-t =n - t$. Thus, $v_j \sim v_\ell$.  This in turn implies that $U^\alpha $ is a clique in $G$, 
since  $U^\alpha \subseteq \{ v_i : d_i \ge \alpha - t, i\in [1,n]\}$. 
\qed

\medskip

\begin{CLA}\label{4.2} Let $\alpha  < \frac{n}{2} $ be any positive integer. If for every $i\in [1,n]$, it holds that $d_i < \alpha - t$ implies $d_i \ge i - t+ 1$,  then $U^\alpha $ is a universal clique in $G$.
\end{CLA}

\proof[Proof of Claim~\ref{4.2}] Assume there exists a positive integer $\alpha < 
\frac{n}{2}$ that satisfies the hypothesis, but  $U^\alpha $ is not a universal clique. Choose $p\in [1,n]$ to be maximum such that there exists $v_q \in U^\alpha$ for some $q\in [1,n]$ such that $v_p\not\sim v_q$. 
By Claim~\ref{4.1}, $v_p \notin \{ v_i : d_i \ge \alpha - t, i\in[1,n]\}$. Thus  $d_p \ge p- t +1$ by the assumption of this claim. By the maximality of $p$, we have $v_q \sim v_\ell$ for all $\ell \in [p+1, n]$. So, $d_q \ge  n - p-1$, which gives $d_p + d_q\ge p- t + 1+n - p - 1  =  n - t$. But, this implies $v_p \sim v_q$, a contradiction. 
\qed

\medskip

Let $\Omega \subseteq V(G)$ be a universal clique in $G$ of maximum size. 

\medskip

\begin{CLA}\label{4.3}  We have $|\Omega| \leq k - 2$.
\end{CLA}

\proof[Proof of Claim~\ref{4.3}] Suppose that $|\Omega| \ge k - 1$.  As $\Omega$ is a universal clique in $G$, we have $d_i \geq | \Omega | \ge k - 1$
for all $i\in [1,n]$. 
 If $| \Omega| > k$, then $d_1 > k$, which contradicts $d_1 \le d_k = k$. Thus  $| \Omega | \le k$.  Note that $v_i\not\in \Omega$ for any $i\in [1,k]$ 
 as every vertex of $\Omega$ has degree $n-1$ and $n-1 > \frac{n}{2}>k$. 
 Let $S=\left(\bigcup_{i\in [1,k]}N(v_i)\right)\setminus \{v_1, \ldots, v_k\}$.   Then we have $\Omega\subseteq S$. 
 % (See Figure~\ref{f0a} for  an illustrations of the relations among $\{v_1,\ldots, v_k\}$, $S$, and $\Omega$.) 
 As  $d_i\le k$ for all $i\in[1,k]$, each $v_i$ has at most $k-|\Omega|$ (this number is either 0 or 1) neighbors  from $\{v_{k+1}, \ldots, v_n\}\setminus \Omega$ in $G$, we have 
 \begin{numcases}{|S| \le } 
 |\Omega| =k &\text{if $|\Omega| = k$}, \nonumber \\
 |\Omega| + k =2k-1  &\text{if $|\Omega| = k-1$}. \nonumber
 \end{numcases}
 Since $\Delta(G[\{v_1, \ldots, v_k\}]) \le 1$, we have 
  $c(G-S) \ge c(G[\{v_1, \ldots, v_k\}]) \ge \frac{k}{2} \ge 4$. However, we get $\frac{|S|}{c(G-S)} <4$, contradicting the toughness of $G$. Thus, Claim~\ref{4.3} must hold.
\qed

\medskip

%\begin{figure}[!htb]
%	\begin{center}
%		\includegraphics[width=0.5\linewidth]{g1.pdf}
%		\caption{The relations among the sets $\{v_1,\ldots, v_k\}$, $S$, and $\Omega$. }
%		\label{f0a}
%	\end{center}	
%\end{figure}

\begin{CLA}\label{4.4}  For all positive integers $\alpha < \frac{n}{2}$ such that $d_\alpha \le \alpha$, we have $|U^\alpha | \ge \alpha - t+1 $.
\end{CLA}

\proof[Proof of Claim~\ref{4.4}] Suppose $v_\alpha \in V(G)$ such that $d_\alpha \le \alpha < \frac{n}{2} $. By the hypothesis of Theorem~\ref{main},  we have $d_{n- \alpha +t} \ge n - \alpha$. That is, there are at least $n - (n - \alpha + t) +1 = \alpha - t +1 $ vertices of degree at least $n - \alpha$, indicating   $|U^\alpha | \ge \alpha - t+1$.
\qed

\medskip

\begin{CLA}\label{4.5} We have $d_\alpha > \alpha$ for all integers  $\alpha$ with  $k+ t - 1 \le \alpha < \frac{n}{2}$.  
\end{CLA}

\proof[Proof of Claim~\ref{4.5}] Assume there exists $\alpha$ such that $k+ t - 1 \le \alpha < \frac{n}{2}$ and $d_\alpha \le \alpha$. Choose such an $\alpha$ to be minimum.  It suffices to show that $U^\alpha$ is a universal clique  in $G$ in order to achieve a contradiction. This is because if   $U^\alpha$ is a universal clique  in $G$, 
then by Claims~\ref{4.3} and ~\ref{4.4},  we have $k - 2 \ge |\Omega| \ge |U^\alpha | \ge \alpha - t $. Rearranging gives $k + t - 2 \ge \alpha  \ge k + t - 1$, a contradiction. 
Thus we show that $U^\alpha$ is a universal clique in $G$  in the following. 
To show that $U^\alpha $ is a  universal clique  in $G$, by Claim~\ref{4.2},  we show that for every $j\in [1,n]$, it holds that $d_j< \alpha - t$ implies $d_j \ge j - t+ 1$. 

We first show  that $d_j\ge \alpha-t$ for all $j\in [\alpha,n]$. 
Consider for now that $j=\alpha$. 
If  $\alpha > k+ t - 1 $, then $\alpha - 1 \ge k+ t - 1$.  By the minimality of $\alpha$,  we get $\alpha - 1 < d_{\alpha - 1} \le d_\alpha \le \alpha$. Thus  $d_\alpha = \alpha >\alpha - t$. If $\alpha = k+t - 1$,  then $d_\alpha \ge d_k = k > \alpha - t $.  In either case, we have shown $d_\alpha > \alpha - t$.  For any $j\in [\alpha+1,n]$, we have $d_j\ge d_\alpha >\alpha-t$. 
Now for  $j\in [1,\alpha-1]$, suppose $d_j< \alpha-t$. 
By the minimality of $k$, we have $d_j  \ge j\ge   j-t+1$  if $j\in [1,k]$.  We have  $d_j \ge d_k = k > k - 1 \ge j - t+1$  if  $j\in [k+1, k+t-2]$. 
By the minimality of $\alpha$, we have $d_j> j>j-t+1$ for all $j\in [k+ t - 1, \alpha-1]$.    This completes the proof. 
\qed

\medskip

\begin{CLA}\label{4.6}  We have $ k \ge \frac{n}{2} - t$.
\end{CLA}

\proof[Proof of Claim~\ref{4.6}]  
We suppose to the contrary that   $k < \frac{n}{2} - t$.  Let  $p = \lfloor \frac{n-1}{2} \rfloor$.  
Then $ k+t-1 \le  p < n/2$.  By Claim~\ref{4.5}, we have $d_p>p$. 
If $d_p = n - 1$, then all vertices from $\{v_p, \ldots, v_n\}$ are contained in a universal clique of $G$ and so we have  $|\Omega | > \frac{n}{2}$. This gives $k  \ge |\Omega | > \frac{n}{2}$, a contradiction  to the assumption that $k < \frac{n}{2} - t$.
Thus there exists  $i\in [1,n]$ such that $v_p\not\sim v_i$.  We choose a maximum such $i$. 
Since $v_i \nsim v_p$,  we have $d_i < n - t - d_p < n - t - (\frac{n-1}{2}-1) = \frac{n+1}{2} - t + 1 \le  d_p$, which gives $i <p$.  
We show that $d_i \ge i - t + 1$. 
If $i \in [1,k]$, then 
by the minimality of $k$, we have $d_i \ge i \ge i - t + 1$. 
If $i \in [k, k + t - 2]$, then 
$d_i \ge d_k = k > k - 1 \ge i - t + 1  $. 
If   $i\in [k+t-1, p-1]$, then  we have $d_i  \ge  i - t + 1$
by Claim~\ref{4.5}. 
Therefore, $d_i \ge i - t + 1$.
By the maximality of $i$, we have $v_p \sim v_j$ for all $j\in [i+1,n]$ and so $d_p \ge n - i - 1$. This gives $d_i + d_p \ge n - i - 1 + i - t + 1 = n - t$, which contradicts that $v_p \nsim v_i$. 
\qed

\medskip

\begin{CLA}\label{4.7}  We have $\delta(G) > \frac{n}{t+1} - 1$. 
\end{CLA}

\proof[Proof of Claim~\ref{4.7}] Assume $\delta(G) \le  \frac{n}{t+1} - 1$. Then, as $2t \le \delta(G)$, we have $(2t+1)(t+1) \le n$. 
By the  minimality of $k$, we have $d_i>i$ for all $i\in [1,k-1]$. For any $i\in [k,n]$, we have $d_i\ge d_k=k >k-t$. 
Thus, by  Claim~\ref{4.2}, we know that $U^k$ is a universal clique in $G$. Therefore, by Claims~\ref{4.4} and~\ref{4.6}, we get  $\delta(G) \ge |U^k| \ge k - t \ge \frac{n}{2} - 2t$.  Observe that for $t \ge 3$,  we have 
\begin{eqnarray*}
	\frac{n}{2}-\frac{n}{t+1} &=& \frac{n(t-1)}{2(t+1)} \ge 
 \frac{(2t+1)(t+1)(t-1)}{2(t+1)}  \\
	 &=&(t+0.5)(t-1) > 2t-1. 
\end{eqnarray*}
This gives  $ \frac{n}{2} - 2t > \frac{n}{t+1} - 1$. Thus $\delta(G) \ge  k -t> \frac{n}{t+1} - 1$, a contradiction. 
\qed

\medskip

As  $\delta(G) > \frac{n}{t+1} - 1$, Theorem~\ref{bauer} implies that $G$ is Hamiltonian, a contradiction to our assumption that $G$ is not Hamiltonian. This completes the proof. 
\qqed

\section{Proof of Theorem~\ref{thm:t-closure}}

\setcounter{THM}{4}

\begin{THM}[Toughness Closure Lemma]\label{thm:t-closure}
	Let $t\ge 4$ be a rational number,   $G$ be  a $t$-tough graph on $n \geq 3$ vertices, and let  distinct   $x, y \in V(G)$  be non-adjacent with $\deg(x) + \deg(y) \geq n - t$. 
	Then $G$ is Hamiltonian if and only if  $G+xy$ is Hamiltonian.  
\end{THM}

\proof It is clear that $G$ being Hamiltonian implies that $G+xy$ is Hamiltonian.  For the converse, we  suppose that $G+xy$ 
is Hamiltonian but $G$ is not.  This implies that $G$ is not complete and so $\delta(G) \ge 2t$. 

As $G+xy$ 
is Hamiltonian, $G$ has a Hamilton path connecting $x$ and $y$. Let $P=v_1\ldots v_n$ be such a path, where $v_1=x$ and $v_n=y$.   We will orient $P$ 
to be from $x$ to $y$, and write $u\preceq v$ for two vertices $u$ and $v$ such that $u$  is at least as close to $x$ along $\oP$ as  $v$ is. 
Our goal is to 
find a cutset $S$ of $G$ with size less than $2t$ and so arrive at a contradiction to the toughness of $G$.  For this purpose, based on the 
assumption that $G$ is not Hamiltonian, we look at how the neighbors of $x$ and $y$ are arranged along this path $P$, and their adjacency relations. 

The first two assertions below follow directly from the assumption that $G$ is not Hamiltonian, and the last two  are corollaries  of the first two. 

\begin{CLA}\label{claim:adjacency}
Let  distinct $i,j \in [2,n-1]$ and suppose $x\sim v_i$ and $y\sim v_j$. Then the following hold. 
\begin{enumerate}[(1)]
	\item  $y\not\sim v_i^-$ and  $v_i^- \not\sim v_j^+$ if  $i<j$. 
	\item If $i>j$, then $v_i^+\not\sim v_j^+$ and $v_i^-\not\sim v_j^-$. 
	\item If $i\le n-3$ and additionally $x\sim v_{i+2}$, then $v_{i+1} \not\sim v_k^+$ for any $v_k$ with $v_k\sim y$. 
	  \item  If $j\le n-3$ and additionally $y\sim v_{j+2}$, then $v_{j+1} \not\sim v_k^-$ for any $v_k$ with $v_k\sim x$. 
\end{enumerate}
\end{CLA}

Since $\deg(x)+\deg(y) \ge n-t$ and $x$ and $y$ do not have two common neighbors that are consecutive on $P$ by Claim~\ref{claim:adjacency}(1) above, 
 each of $x$ and $y$ is expected to have many neighbors that are consecutive on $P$.  Thus we define neighbor intervals for $x$ and $y$, respectively, as a set of consecutive vertices on $P$ that are all adjacent to $x$ or $y$.  
  For $z\in \{x,y\}$, and $v_i, v_j$ with $i,j\in [2,n-1]$  and $i\le j$  such that $z\sim v_i, v_j$,  we call $V(v_iPv_j)$ a \emph{$z$-interval} and denote it by 
  $I_z[v_i, v_j]$ if $V(v_iPv_j) \subseteq N(z)$ but  $v_i^-, v_j^+ \not\sim z$.   
%  We let  $$I_z(v_i, v_j) =I_z[v_i, v_j]\setminus \{v_i, v_j\}, \quad I_z(v_i, v_j] =I_z[v_i, v_j]\setminus \{v_i\}, 
% \quad I_z[v_i, v_j) =I_z[v_i, v_j]\setminus \{ v_j\}.$$
 
Given $I_x[v_i,v_j]$ and $I_y[v_k,v_\ell]$, by Claim~\ref{claim:adjacency}(1), we know that the two intervals can have at most one vertex in common. In case that they do have 
a common vertex, then it must be the case that $v_j=v_k$.  In this case, we let $I_{xy}[v_i,v_j,v_\ell]=I_x[v_i,v_j]\cup I_y[v_k,v_\ell]$ and call it a \emph{joint-interval}. 
For $i,j\in [3,n-2]$ with $i\le j$,  we define an \emph{interval-gap} to be a set of consecutive vertices on $P$ that are all adjacent to neither  $x$ nor $y$. 
A  \emph{parallel-gap} is $J[v_i,v_j]:= V(v_iPv_j)$ such that $V(v_iPv_j)\cap (N(x)\cup N(y))=\emptyset$ and  that $v_i^-, v_j^+ \in N(x)$,  or $v_i^-, v_j^+ \in N(y)$, 
or $v_i^-\in N(x)$ but $v_j^+\in N(y)$. A   \emph{crossing-gap} is  $J[v_i,v_j]:= V(v_iPv_j)$ such that $V(v_iPv_j)\cap (N(x)\cup N(y))=\emptyset$ and  that $v_i^-\in N(y)$ and $v_j^+\in N(x)$. 
By the range of  $i$ and $j$ in the above definition, we see that each of $x$ and $y$  is not contained in any interval-gaps.

Let $\mathcal{I}_x$ be the set of $x$-intervals that are not joint-intervals, $\mathcal{I}_y$ be the set of $y$-intervals that are not joint-intervals, and $\mathcal{I}_{xy}$ be the set of joint-intervals. 
Let 
$$
p=|\mathcal{I}_x \cup \mathcal{I}_y|, \quad \text{and}\quad q= |\mathcal{I}_{xy}|. 
$$

\begin{CLA}\label{claim:crossing-gap}
Each crossing-gap contains at least two vertices and there are  at least $q-1$ distinct crossing-gaps when $q\ge 1$. 
\end{CLA}

\proof[Proof of Claim~\ref{claim:crossing-gap}]   For the first part, suppose $\{v_i\}$ for some $i\in [2,n-1]$ is a crossing-gap with a single vertex.  Then $C=v_{i+1}x\oP v_{i-1} y \iP v_{i+1}$ gives a Hamilton cycle of $G-v_i$. 
We have $v_i\sim v_{i-1}, v_{i+1}$, and 
with respect to the cycle $\oC$, we have $x=v_{i+1}^+$ and $y=v_{i-1}^+$. However, $\deg(x)+\deg(y) \ge n-t$, contradicting Theorem~\ref{thm:sum-of-successors}. 
For the second part, assume that $q\ge 2$. Let the $q$ common neighbors of $x$ and $y$ be $u_1, \ldots u_q$ with $u_1\preceq u_2 \ldots \preceq u_q$.  Thus 
$V(u_iPu_{i+1})$ for  each $i\in [1,q-1]$ is a set of vertices such that $u_i\sim y$ and $u_{i+1}\sim x$. 
We let $w_1 \in V(u_iPu^-_{i+1})$ with $w_1$ farthest away from $u_i$  along $\oP$ such that  $w_1 \sim y$ (note that $w_1$ exists as $u_i\sim y$). 
By the choice of $w_1$, we know that $V(w_1^+P u^-_{i+1}) \cap N(y) =\emptyset$. 
Then we let $w_2 \in V(w^+_1Pu_{i+1})$ with $w_2$ closest to  $w^+_1$  along $\oP$ such that  $w_2 \sim x$ (note that $w_2$ exists as $u_{i+1}\sim x$). 
By the choice of $w_2$, we know that $V(w_1^+P w_2^-) \cap N(x) =\emptyset$.  
Therefore, $V(w^+_1Pw^-_2)\cap (N(x)\cup N(y)) =\emptyset$.  Note that $w_1\ne w_2$ as $w_2 \in V(w^+_1Pu_{i+1})$. 
 Thus  $V(w^+_1Pw^-_2)$ is a crossing-gap. 
Since $V(u_i^+Pu^-_{i+1})$ and $V(u_j^+Pu^-_{j+1})$ are disjoint for distinct $i, j\in [1,q-1]$, 
we can find $q-1$ distinct crossing-gaps. 
\qed 

Let $p^*$ be the total number of distinct parallel-gaps and $q^*$ be the total number of distinct crossing-gaps.   We  let  the  set of $p^*$ parallel-gaps be $\{J[u_i,w_i]: i\in [1,p^*], u_1\preceq w_1 \preceq u_2 \preceq w_2 \preceq \ldots \preceq u_{p^*} \preceq w_{p^*}\}$, and let  $|J[u_i,w_i]|=p_i$. We also let  the  set of $q^*$ crossing-gaps be $\{J[r_i,s_i]: i\in [1,q^*], r_1\preceq s_1 \preceq r_2 \preceq s_2 \ldots \preceq r_{q^*} \preceq s_{q^*}\}$, and let $|J[r_i,s_i]|=q_i$.

\begin{CLA}\label{2.2.1}   We have $|\mathcal{I}_x \cup \mathcal{I}_y \cup  \mathcal{I}_{xy}|=p+q\leq t - \sum\limits_{i=1}^{p^*} (p_i -1) - \sum\limits_{i =1 }^{q^*} (q_i -2)$. 
\end{CLA}

%Case q=0 

\proof[Proof of Claim~\ref{2.2.1}]    By the definition,  the three sets $\mathcal{I}_x, \mathcal{I}_y, \mathcal{I}_{xy}$ are pairwise disjoint. Thus 
$|\mathcal{I}_x \cup \mathcal{I}_y \cup  \mathcal{I}_{xy}| =p+q$.  Also, by our definition, we have $|N(x)\cap N(y)| =|\mathcal{I}_{xy}|=q$ and so $|N(x) \cup N(y)|  \ge  n - t  - q$.
Since $|\mathcal{I}_x \cup \mathcal{I}_y \cup  \mathcal{I}_{xy}| =p+q$, and $v_2$ and $v_{n-1}$ are contained in an $x$-interval, $y$-interval, or  joint-interval, it follows that  there are exactly $p+q-1 =p^*+q^*$ interval-gaps. By Claim~\ref{claim:crossing-gap},  $q^*\ge q-1$. 
As each of $x$ and $y$ is not contained in any interval-gaps, and is not contained in $N(x)\cup N(y)$ by the assumption that $x\not\sim y$, we get 
\begin{eqnarray*}
	t+q &\ge & |V(G) \setminus  (N(x) \cup N(y))|   \ge  2 +\sum_{i=1}^{p^*} p_i  + \sum_{i = 1}^{q^*}q_i \\
	 &\ge & 2 + p^* + \sum\limits_{i=1}^{p^*} (p_i-1)  + 2q^*  + \sum\limits_{i = 1}^{q^*}(q_i-2). 
\end{eqnarray*}
As $p+q-1 =p^*+q^*$ and $q^*\ge q-1$, we get $p+q \le t  - \sum\limits_{i=1}^{p^*} (p_i-1)   - \sum\limits_{i = 1}^{q^*}(q_i-2)$. 
Therefore, 
\[
|\mathcal{I}_x \cup \mathcal{I}_y \cup  \mathcal{I}_{xy}|=p+q \leq t  - \sum_{i=1}^p (p_i-1)   - \sum_{i = 1}^{q-1}(q_i-2), 
\]
as desired. 
\qed

\begin{CLA}\label{2.2.2}  For any $i\in[2,n-2]$, if $\{v_i,v_{i+1}\}$ is a crossing-gap of size 2, then $v_i \not\sim v_j$
	for any $j\in[3,n-2]$ such that $y\sim v_{j-1}, v_{j+1}$. 
\end{CLA}

\proof[Proof of Claim~\ref{2.2.2}] Suppose to the contrary that $v_i\sim v_j$. We will show that $v_{i+1}$ has less than $2t$ neighbors in $G$ to arrive at a contradiction to $G$ being $t$-tough.

\begin{figure}[!htb]
	\begin{center}
		\includegraphics[width=1\linewidth]{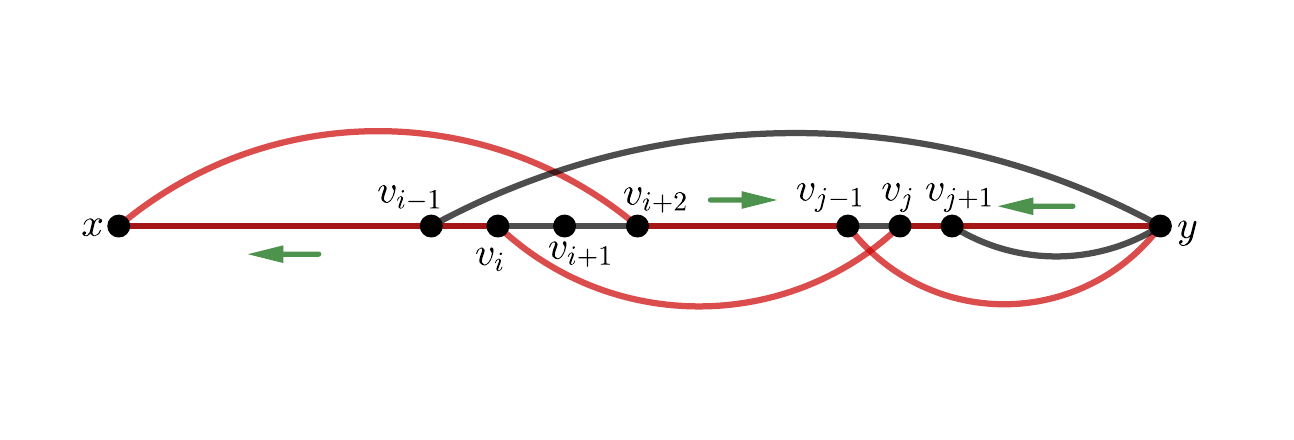}
		\caption{Construction of $C$ when $i<j$, drawn in red. The green arrows indicate the orientation of the corresponding segments of $P$ on $\oC$. }
		\label{f0}
	\end{center}	
\end{figure}

By Claim~\ref{claim:adjacency}(1)-(2), we know that for any $v_k\sim y$ with $v_k\preceq v_i$ on $P$, we have $v_{i+1}\not\sim v_{k-1}$;
and for  $v_k\sim y$ with $ v_i \preceq v_k$ on $P$, we have $v_{i+1}\not\sim v_{k+1}$.  
Thus   vertices from $(N(y)\cap V(v_2Pv_i))^-$ and $(N(y)\cap V(v_{i+2}Pv_{n-1}))^+$
are non-neighbors of $v_{i+1}$. 
Let 
\begin{numcases}{C=}
v_j v_i \lP x v_{i+2} \rP v_{j-1} y \lP v_j & \text{if $i<j$ (see Figure~\ref{f0})}, \nonumber \\
v_j v_i \lP v_{j+1} y \lP v_{i+2} x \rP v_j & \text{if $i>j$}. \nonumber
\end{numcases}
Then $C$ is a Hamilton cycle of $G-v_{i+1}$.    The predecessors and successors 
of vertices below are all taken with respect to $\oC$. 
As $G$ is not Hamiltonian,  both $N(v_{i+1})^-$ and $N(v_{i+1})^+$ are independent in $G$. 
When $i<j$, since $v_{i+1}\sim v_{i+2}$ and $x=v_{i+2}^-$, it then follows that $v_{i+1}\not\sim 
z^+$ for any $z\in N(x)$. As a consequence, we get $N(x)^+ \cap N(v_{i+1}) = \emptyset$. 
When $i>j$, since $v_{i+1}\sim v_{i+2}$ and $x=v_{i+2}^+$, it then follows that $v_{i+1}\not\sim 
z^-$   for any $z\in N(x)$. As a consequence, we get $N(x)^- \cap N(v_{i+1}) = \emptyset$. 

%The arguments above indicate that for every distinct vertex $z\in N(x)\cup N(y)$, 
%there is a unique non-neighbor of $v_{i+1}$ that is corresponding to $z$. 
%Thus $v_{i+1}$ has at least 
%$|N(x)\cup N(y)|$  non-neighbors on $C$.  Then  by Claim~\ref{2.2.1} that $q\le t$, we get 
%\begin{eqnarray*}
%	\deg(v_{i+1}) &\le & n-1-|N(x)\cup N(y)|  \\
%	& \le&  n-1-(n-t-q)   \\
%	& \le&  2t-1, 
%\end{eqnarray*}
%%% The line above is false as a vertex from $N(x)$ and a vertex from $N(y)$ can be corresponding to the same non-neighbor of $v_{i+1}$. Thus we need a different argument. 
When $i<j$, by the construction of $C$ and the arguments before, we have  $v_{i+1}\not\sim z^+$ for any $z\in (N(x)\cup N(y))\cap (V(xPv_i) \cup V(v_{i+2}Pv_{j-1}))$,   $v_{i+1}\not\sim z^+$ for $z\in N(x)\cap V(v_jPy)$, 
 $v_{i+1}\not\sim z^-$ for $z\in N(y)\cap V(v_jPy)$, any joint-interval  contained in $V(v_jPy)$  has the corresponding 
 $x$-interval preceding the corresponding 
 $y$-interval along $\oC$, 
and  there is no joint-interval with some vertices  in 
$V(xPv_i) \cup V(v_{i+2}Pv_{j-1})$ and some other vertices  in $V(v_jPy)$. Thus 
 $v_{i+1}$ can have at most one neighbor from each set in $\mathcal{I}_x \cup \mathcal{I}_y \cup  \mathcal{I}_{xy}$,
 which holds also true when $i>j$ 
 by following the same argument. 
 
For each interval-gap, say $\{w\}$,  of size one,  we claim that we can assume $v_{i+1}\not\sim w$. 
 We only consider the case $i<j$, as the argument for the other case follows the same logic.  By the construction of $C$, if $w^-, w^+\in V(xPv_i) \cup V(v_{i+2}Pv_{j-1})$, then 
 $v_{i+1}\not\sim w$ since  $w^-\in (N(x)\cup N(y))\cap (V(xPv_i) \cup V(v_{i+2}Pv_{j-1}))$.

 Consider then  that $w^- \in V(v_{i+2}Pv_{j-1})$ and $w^+ \in V(v_jPy)$. 
 Then we have $w^-=v_{j-1}$ and $w^+=v_{n-1}$, and $w=y$. As $v_{i+1}$ is a vertex from an interval-gap, 
 we have $v_{i+1} \not\sim w$.

 Lastly consider  $w^-, w^+ \in V(v_jPy)$.  Then we have $v_{i+1}\not\sim z^+$ for $z\in N(x)\cap V(v_jPy)$ 
 and $v_{i+1}\not\sim z^-$ for $z\in N(y)\cap V(v_jPy)$. 
Thus we have $v_{i+1}\not\sim w$ if $x\sim w^-$. Hence  we assume that $y\sim w^-$. If $y\sim w^+$, then $v_{i+1}\not\sim w$. Thus we assume that $x\sim w^+$. 
 This implies that  $w^+$ is the only possible neighbor of  $v_{i+1}$ from vertices in the $x$-interval containing  $w^+$ and $w^-$ is the only possible neighbor of  $v_{i+1}$ from vertices in the $y$-interval containing  $w^-$. If $v_{i+1} \sim w^-$ or $v_{i+1} \sim w^+$, then $v_{i+1}\not\sim w$ as $v_{i+1}$ 
 has no two consecutive neighbors on $C$.  Thus,  we assume that $v_{i+1} \not\sim w^-$ and $v_{i+1} \not\sim w^+$. 
 This implies that $v_{i+1}$ has no neighbor in the  $x$-interval containing  $w^+$, and has no neighbor in  the $y$-interval containing  $w^-$.  
 Therefore, in our counting of $\deg(v_{i+1})$,  if it does hold that $v_{i+1} \sim w$, we   distribute  the contribution of $w$ 
  in $\deg(v_{i+1})$  by assuming that  $v_{i+1}$  has a neighbor in the  $x$-interval containing  $w^+$.  Thus, 
from this sense of counting, we may assume $v_{i+1}\not\sim w$. 
 
As $v_{i+1}$ 
has no two consecutive neighbors on $C$, the above arguments indicate that $v_{i+1}$ has at most  $\frac{1}{2}(n-1-|N(x)\cup N(y)|)$ neighbors from $V(G-v_{i+1})\setminus (N(x)\cup N(y))$. 
As $|\mathcal{I}_x \cup \mathcal{I}_y \cup  \mathcal{I}_{xy}| =p+q \le t$ by Claim~\ref{2.2.1} and $|N(x)\cup N(y)| \ge n-( t+q)$, we know that 
\begin{eqnarray*}
	\deg(v_{i+1}) &\le & |\mathcal{I}_x \cup \mathcal{I}_y \cup  \mathcal{I}_{xy}| +\frac{1}{2}\big(n-1-|N(x)\cup N(y)| \big) \\
	& \le&  t+ \frac{1}{2}(t+q-1)  \\
	& <&  2t, 
\end{eqnarray*}
a contradiction. 
\qed

We now construct a cutset $S$  of $G$ such that $|S| < 2t$. To do so, we define the following sets:
\begin{eqnarray*}
S_x&=&\{v_j, v_{j+1}:  \text{$v_j$ is the right endvertex of an $x$-interval that is not a joint-interval}\},  \\
S_y&=&\{v_i, v_j:  \text{$I_y[v_i,v_j]$ is   a $y$-interval that  is not a joint-interval}\}, \\
S_{xy} &=& \{v_j, v_\ell: \text{$I_{xy}[v_i,v_j,v_\ell]$ is a joint-interval}\}, \\
T_1&=& \bigcup_{\text{ $J[v_i,v_j]$ is a parallel-gap of size at least 2}} J[v_i,v_j], \\ 
T_2&=& \bigcup_{\text{ $J[v_i,v_j]$ is a crossing-gap of size 3}} \left(J[v_i,v_j]\setminus\{v_j\}\right),  \\
T_3&=& \bigcup_{\text{ $J[v_i,v_j]$ is a crossing-gap of size at least 4}} J[v_i,v_j]. 
\end{eqnarray*}

Let 
\begin{numcases}{S=}
S_x\cup S_y\cup S_{xy}\cup T_1\cup T_2 \cup  T_3 & \text{if  $\{v_{n-1}\}$ is a $y$-interval}, \nonumber \\ 
\left(S_x\cup S_y\cup S_{xy}\cup T_1\cup T_2 \cup  T_3\right)\setminus \{v_{n-1}\} & \text{otherwise}. \nonumber 
\end{numcases}
We 
prove the following claims regarding what vertices are in $V(G)\setminus S$ and the size of $S$. 

\begin{CLA}\label{claim:left-over}
 Let $v_i\in V(G)\setminus S$ for some $i\in [2,n-2]$. Then $x\sim v_i, v_{i+1}$, or $y\sim v_{i-1}, v_{i+1}$, 
 or $v_i$  is contained in a parallel-gap of size one such that $y\sim v_{i-1}, v_{i+1}$, 
 or $v_i$ is contained in a crossing-gap of size two, or $v_i$ is the right  endvertex of a crossing-gap of size three. 
\end{CLA}

\proof[Proof of Claim~\ref{claim:left-over}]  By the definition of $S$, we know that either $v_i$ is a neighbor of $x$ or $y$, or 
 $v_i$ is contained in 
a parallel-gap of size one, or a crossing-gap of size two or three.  If $x\sim v_i$, then by the definition of $S_x$, we have $x\sim v_{i+1}$. 
If $y\sim v_i$, then by the definition of $S_y$, we have $y\sim v_{i-1}, v_{i+1}$. 
If $v_i$ is contained in a parallel-gap of size one, then by the definition of $S_x$, we know that $y\sim v_{i-1}$. As $\{v_i\}$ is a parallel-gap, 
$y\sim v_{i-1}$ implies  $y\sim v_{i+1}$. If $v_i$ is contained in crossing-gap of size   three, 
then $v_i$ is the right endvertex of a crossing-gap of size three by the definition of $T_3$. 
\qed 

\begin{CLA}\label{claim:S-size}
	We have $|S| \le 2t-1$. 
\end{CLA}

\proof[Proof of Claim~\ref{claim:S-size}] 

  For each crossing-gap $J[r_i,s_i]$ of size $q_i$, we let $q_i^* =q_i$ if $q_i\ge 4$, $q_i^* =q_i-1$ if $q_i= 3$, 
and $q_i^* =0$ if $q_i= 2$.  Note that by the definition of $S$, only one vertex was deleted from the $y$-interval containing $v_{n-1}$.
Now  by the definition of $S$ and Claim~\ref{2.2.1}, we have 
\begin{eqnarray*}
	|S| &\le & 2(p+q) -1 + \sum_{i=1, p_i\ge 2}^{p^*} p_i + \sum_{i=1}^{q^*} q_i^*   \\ 
	 &\le &  2\left(t - \sum\limits_{i=1}^{p^*} (p_i -1) - \sum\limits_{i =1 }^{q^*} (q_i -2)\right)-1+\sum_{i=1, p_i\ge 2}^{p^*} p_i + \sum_{i=1}^{q^*} q_i^*   \\
	 &=& 2t-1+\sum_{i=1, p_i\ge 2}^{p^*} \left(p_i-2(p_i-1)\right)+\sum\limits_{i =1 }^{q^*}\left(q_i^*-2(q_i-2)\right)\\
	 &\le & 2t-1, 
\end{eqnarray*}
where the last inequality follows as $p_i-2(p_i-1)  \le 0$ when $p_i\ge 2$, and $q_i^*-2(q_i-2) \le 0$ by the definition of $q_i^*$
and the fact that $q_i\ge 2$ for all $i\in [1,q^*]$ from Claim~\ref{claim:crossing-gap}. 
\qed

\begin{CLA}\label{claim:S-cut}
	We have $c(G-S) \ge 2$. 
\end{CLA}

\proof[Proof of Claim~\ref{claim:S-cut}]
For the sake of contradiction, suppose $G' = G - S$ is connected. Let $X' = N_{G'}(x) \cup \{x\}$ and $Y' = N_{G'}(y) \cup \{y\}$. Then, there must exist
a path  $P'$ in $G'$ connecting a vertex of $X'$ and a vertex $Y'$ which is internally-disjoint with $X'\cup Y'$. Suppose that $P'=uu_1\ldots u_h v$ 
for some $u\in X'$ and $v\in Y'$.  By Claim~\ref{claim:left-over},  we know that $v=y$,  or $V^{\leftharpoonup}, v^+ \sim y$, 
or $v=v_{n-1}$ when the $y$-interval containing $v_{n-1}$ has size at least two,   and that $u^+\sim x$. 
By Claim~\ref{claim:adjacency}(1) and (4), we know that 
$P'\ne uv$. Thus $P'$ contains at least three vertices.  As $P'$ is internally-disjoint with $X'\cup Y'$, $u_1, \ldots, u_h$
are from interval-gaps of $P$.  

Recall again,  $v=y$,  or $v^-, v^+ \sim y$, 
or $v=v_{n-1}$ when the $y$-interval containing $v_{n-1}$ has size at least two. Since $u_h\sim v$,    Claim~\ref{claim:adjacency}(1) and (4) imply  that $u_h^+ \not\sim x$. 
Thus $u_h$ is not the right endvertex of any crossing-gap. By Claim~\ref{2.2.2}, $u_{h}$ is not the 
left endvertex of any crossing-gap of size two. Thus by Claim~\ref{claim:left-over}, $\{u_h\}$
is a parallel-gap of size one such that $y\sim u_h^-, u_h^+$. Now with $u_h$ in the place of $v$,  the same arguments 
as above imply that $\{u_{h-1}\}$, if it exists,  is a parallel-gap of size one such that $y\sim u_{h-1}^-, u_{h-1}^+$. 
Similarly, for any $i\in [1,h-2]$,   if it exists, we deduce that $\{u_{i}\}$  is a parallel-gap of size one such that $y\sim u_{i}^-, u_{i}^+$. 
As $u_1\sim u$ and $u^+\sim x$, we get a contradiction to Claim~\ref{claim:adjacency}(4). 
\qed 

Now Claims~\ref{claim:S-size} and~\ref{claim:S-cut} together give a controduction to the toughness of $G$, 
completing the proof of Theorem~\ref{thm:t-closure}. 
\qqed

\section{Proof of Theorem~\ref{thm:sum-of-successors}}

\setcounter{THM}{5}

\begin{THM}\label{thm:sum-of-successors}
	Let $t\ge 4$ be rational and  $G$ be  a $t$-tough graph on $n \geq 3$ vertices. Suppose that $G$ is not Hamiltonian, but there exists $z \in V(G)$ such that $G -z$ has a Hamilton cycle $C$. Then, for any distinct $x, y \in (N(z))^+$, we have that $\deg(x)+\deg(y) < n-t$. 
\end{THM}

\proof  Suppose to the contrary that  there are distinct $x, y \in (N(z))^+$ for which  $\deg(x)+\deg(y)  \ge  n-t$.  As $G$ is not Hamiltonian, $G$ is not a complete graph. 
	Thus $\deg(z) =\deg(z, C)\geq 2t$.    
	As $C$ is a Hamilton cycle of $G-z$, $G$ has  a Hamilton $(x,y)$-path $P:=x\oC y^-zx^-\lC y$, as demonstrated in 
	Figure~\ref{f3} below.

	\begin{figure}[!htb]
		\begin{center}
			\includegraphics[width=1.1\linewidth]{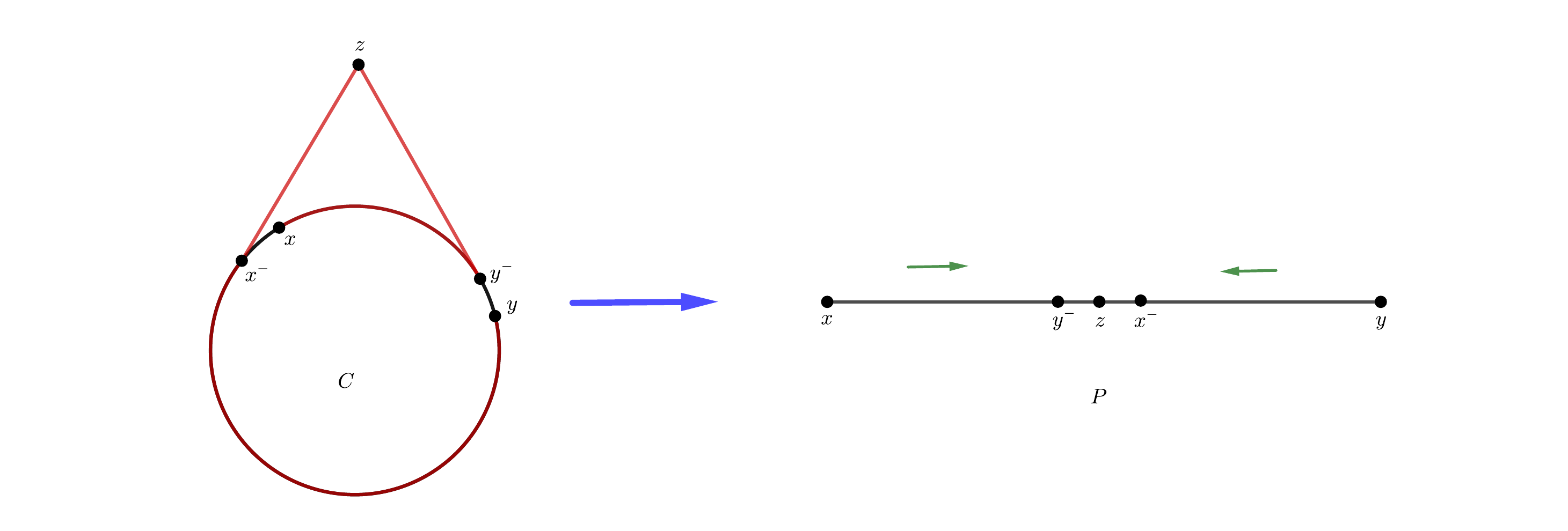}
			\caption{A Hamilton $(x,y)$-path  $P$ of $G$ constructed from $C$, where the first green arrow indicates that the direction from $x$ to $y^-$  on $P$ consists with the clockwise direction of $C$, and the second green arrow indicates that the direction from $x^-$ to $y$  on $P$  is opposite to the clockwise direction of $C$. }
			\label{f3}
		\end{center}	
	\end{figure}

For	 $W, T \subseteq V(G)$, we let $N(W,T)=(\bigcup_{v\in W} N(v) ) \cap W$. 
If $H\subseteq V(G)$, we write $N(W, H)$ for $N(W, V(H))$. 
We fix the forward direction of $P$ to be from $x$ to $y$. For $v \in V(P)$, we let $v^{+}$ denote the neighbor of $v$ on $P$ that precedes $v$ in the forward direction of $P$, and let $v^{-}$ denote the neighbor of $v$ on $P$ that follows $v$ along the direction of $P$. For $S\subseteq V(P)$,  let $S^+=\{v^+: v\in S\}$ and $S^-=\{v^-: v\in S\}$. 
Note that the superscripts of the two special vertices $x^{-}$ and $y^{-}$ are defined with respect to the clockwise direction of $C$.

    We will construct a cutset $S$ of $G$ such that $\frac{|S|}{c(G-S)} <t$. For this purpose, 
    we   define the following sets: 
    \begin{eqnarray*}
    	X&=& \{v\in V(P): v\sim x\}, \\
    	 Y&=&\{v^+\in V(P): v\sim y\}, \\
    	Z&=& \{v^+\in V(xpy^-): v\sim z\} \cup \{v^-\in V(x^-Py): v\sim z\} \cup \{x,y\}, \\
    	R&=& V(G)\setminus (X\cup Y \cup Z). 
    \end{eqnarray*}    
In the following, we prove some properties of these sets. 
    
    \medskip 
    
      \begin{CLA}\label{2.1.1} We have $ X \cap Y = \emptyset$.
      \end{CLA}
\proof[Proof of Claim~\ref{2.1.1}] Suppose to the contrary that there exists $v\in X\cap Y$. Then 
$xvPyv^-Px$ is a Hamilton cycle. 
\qed 

Note that the set $Z$ is precisely the set of the predecessors of the neighbors of $z$ on $\oC$. 
If there are  $u,v \in Z$ such that $u\sim v$, then  $u v \rC u^- z v^- \lC u$ is a Hamilton cycle in $G$.  
Thus we have the following claim. 

\medskip 

\begin{CLA}\label{2.1.2}  The set  $Z$ is an independent set in $G$. 
\end{CLA}
   \medskip

\begin{CLA}\label{2.1.3}  We have  $|R \cup (Z \setminus  Y) | \leq t$ and $|Y\cap Z| \geq |R| +t$. 
\end{CLA}

    \proof[Proof of Claim~\ref{2.1.3}] 
        Clearly $|X\cup Y \cup Z | \leq n - |R|$.  Note that 
         $|X| = \deg(x)$ and $|Y| =\deg(y)$ as $y\not\sim y$. 
          By Claim~\ref{2.1.1}, we have $|X \cup Y| = |X| + |Y|$; and by Claim~\ref{2.1.2}, we have $X \cap Z = \emptyset$.  Thus  we get 
        \begin{eqnarray}
        	  n-|R|\ge |X\cup Y \cup Z | & \ge & |X| + |Y| + |Z| - |X \cap Z| - |Y \cap Z|  \nonumber  \\
        	  &\ge & n-t + |Z| - |Y \cap Z| = n-t+|Z\setminus Y|,  \label{eqn:eqn2}
        \end{eqnarray}
     which gives $|R \cup (Z\setminus Y) | \leq t$.
     For the second part, it follows from~\eqref{eqn:eqn2} by noting that $|Z| \ge 2t$. 
  \qed 
  
Let 
$$
U=\{v\in V(P): v^+\sim x, v^-\sim y\} =X^-\cap Y. 
$$
In the remainder of the proof, we  show that  $U$ is an independent 
set of size at least $t-1$,   $|N(U)|  \le 2t+2|U|$, and $c(G-N(U)) \ge |U|+1$. These properties together contradict the assumption that  $G$ is 4-tough.

\begin{CLA}\label{2.1.4a}  We have  $y\not\in U$ and $y\not\in N(U)$.
\end{CLA}
\proof[Proof of Claim~\ref{2.1.4a}]  
Since $y\not\in X^-$, we have  $y\not\in U$.  As $U=X^-\cap Y$, Claim~\ref{2.1.1}
implies that $y\not\in N(U)$. 
\qed

\begin{CLA}\label{2.1.4}  We have  $|U| \geq t-1$.
 \end{CLA}
 
 \proof[Proof of Claim~\ref{2.1.4}] 
 As $|Y \cap Z| \geq |R| +t$, it suffices to show  $|(Y \cap Z) \setminus  U| \leq |R|+1$.  We show that except possibly the vertex $y$, for every 
 $v\in (Y \cap Z) \setminus  U$  with $v\ne y$,  we have $v^+\in R$. 
 
 Consider first that $v\in V(x^+Py)$. 
 Then as $v\in Y\cap Z$, we have $v^-\sim y$ and  $v^+\sim z$. Since $v\not\in U$, we have $v^+\not\sim x$.
 Thus   $v^+\not\in X$. Since $Z$ is an independent set in $G$ by Claim~\ref{2.1.2} and $v,y\in Z$, 
 it follows that $y\not\sim v$  and $z\not\sim v^{++}$. Thus $v^+\not\in Y\cup Z$. Since  $v^+\not\in X$ also, 
 we have $v^+\in R$.

Consider then  that $v\in V(xPy^-)$. Note that $v\ne x$ as $x \not\in Y$.   Then $v^-\sim y,z$.  
Since $v\not\in U$, we have $v^+\not\in X$. 
Since $Z$ is an independent set in $G$ by Claim~\ref{2.1.2} and $v,y\in Z$, 
it follows that $y\not\sim v$  and $z\not\sim v$. Thus $v^+\not\in Y\cup Z$. Since  $v^+\not\in X$ also, 
we have $v^+\in R$. 
 Therefore we have 
 \begin{eqnarray*}
 |(Z \cap Y) \setminus  U|  & \le & 1+ \left| \Big(\big((Y \cap Z) \setminus \{y\}\big)\setminus U \Big)^+\right| \\
 & \le & 1+|R|. 
 \end{eqnarray*}
 \qed

%\medskip 
%\begin{CLA}\label{2.1.5} The set   $U\cup \{z\}$ is an independent set in $G$.
%\end{CLA}
%
% \proof[Proof of Claim~\ref{2.1.5}]   Note that $u\ne y$ for any $u\in U$ by Claim~\ref{2.1.4a}. 
% Since $x,y\in Z$, $U=  X^-\cap Y$, and $Z$ is an independent set by Claim~\ref{2.1.2}, 
% it follows that $z$ is not adjacent in $G$ to any vertex of $U$ by the definition of $Z$. 
% Next, let distinct $u, v\in U$ such that $u\sim v$.  Suppose, without loss of generality, that 
% $u$ is on $xPv$.  Then $uvPyv^-Pu^+xPu$ is a Hamilton cycle of $G$. 
%Therefore, $U\cup \{z\}$ is an independent set in $G$.
% \qed 
% 	

\begin{CLA}\label{2.1.7}  
	We have  $|N(U ) | \le 2t+2|U|$.  
\end{CLA}
\proof[Proof of Claim~\ref{2.1.7}]

Since $|N(U,P)| \le 2|U|$,  it suffices to show that 
$|N(U)\setminus N(U,P)| \le 2t$.   Since $|(Z\setminus Y) \cup R| \le t$ 
by  Claim~\ref{2.1.3}, it remains to show 
 that every vertex $v\in N(U)\setminus N(U,P)$ corresponds to 
a unique vertex  $g(v)$ of $(Z\setminus Y) \cup R$ and that for every vertex $w\in (Z\setminus Y) \cup R$, 
there are at most two distinct vertices $v,v'\in N(U)\setminus N(U,P)$ such that $g(v)=g(v')=w$. 
 
 Let $v\in N(U)\setminus N(U,P)$ and $u\in U$ such that $u\sim v$.  We consider two cases regarding the locations
 of $u$ and $v$ on $P$ for showing that every vertex $v\in N(U)\setminus N(U,P)$ corresponds to 
 a unique vertex  $g(v)$ of $(Z\setminus Y) \cup R$.

{\bf \noindent Case 1: $u$ is on $xPv$.}

Then we have $v^-\not\sim x$. For otherwise, $xv^-PuvPyu^-Px$  is a Hamilton cycle.   This implies that $v^-\not\in X$. 
Similarly, we have  $v^-\not\sim y$. For otherwise, $yv^-Pu^+xPuvPy$  is a Hamilton cycle.   

If $v^-\not\in Y$, then we have $v^-\in (Z\setminus Y) \cup R$. In this case, we let $g(v)=v^-$. 
Thus we assume that $v^-\in Y$ and so $y\sim v^{--}$.  Then we have $x\not\sim v$ as $v\not\in N(U,P)$ 
implies that $v^-\not\in U$. Thus $v\not\in X$.  Since $v^-\not\sim y$  implies that $v\not\in Y$, 
we then know that $v\in (Z\setminus Y) \cup R$. In this case, we let $g(v)=v$.

{\bf \noindent Case 2: $v$ is on $xPu$.}

Then we have $v^+\not\sim x$. For otherwise, $xv^+Pu^-yPuvPx$  is a Hamilton cycle.   This implies that $v^+\not\in X$. 
Similarly, we have  $v^+\not\sim y$. For otherwise, $yv^+PuvPxu^+Py$  is a Hamilton cycle.   

If $v^+\not\in Y$, then we have $v^+\in (Z\setminus Y) \cup R$. In this case, we let $g(v)=v^+$. 
Thus we assume that $v^+\in Y$ and so $y\sim v$.  Then we have $x\not\sim v^{++}$ as $v\not\in N(U,P)$ 
implies that $v^+\not\in U$. Thus $v^{++}\not\in X$.  Since $v^+\not\sim y$  implies that $v^{++}\not\in Y$, 
we then know that $v^{++}\in (Z\setminus Y) \cup R$. In this case, we let $g(v)=v^{++}$.

We lastly verify that for any $w\in (Z\setminus Y) \cup R$, there are  at mos 
two distinct vertices $v_1,v_2\in N(U)\setminus N(U,P)$ such that $g(v_1)=g(v_2)=w$. 
Suppose to the contrary that  there exists $w\in (Z\setminus Y) \cup R$ 
such that there are three distinct vertices $v_1, v_2, v_3\in N(U)\setminus N(U,P)$ for which $g(v_1)=g(v_2)=g(v_3)=w$.    Let $u_1, u_2, u_3 \in U$ such that $u_i\sim v_i$ for each $i\in [1,3]$. 
Note that $u_1, u_2, u_3$ do not need to be all distinct. 
Then by the Pigeonhole Principle,  without loss of generality, 
 either $u_1, u_2$ is on $xPw$ or $u_1,u_2$ is  on $wPy$, where $w=g(v_1)=g(v_2)$. 
By symmetry, we suppose that $u_1$ is on $xPu_2$.    We will find a Hamilton cycle in each case
and so get a contradiction to the assumption that $G$ is not Hamiltonian. 

\begin{figure}[!htb]
	\begin{center}
		\includegraphics[width=1\linewidth]{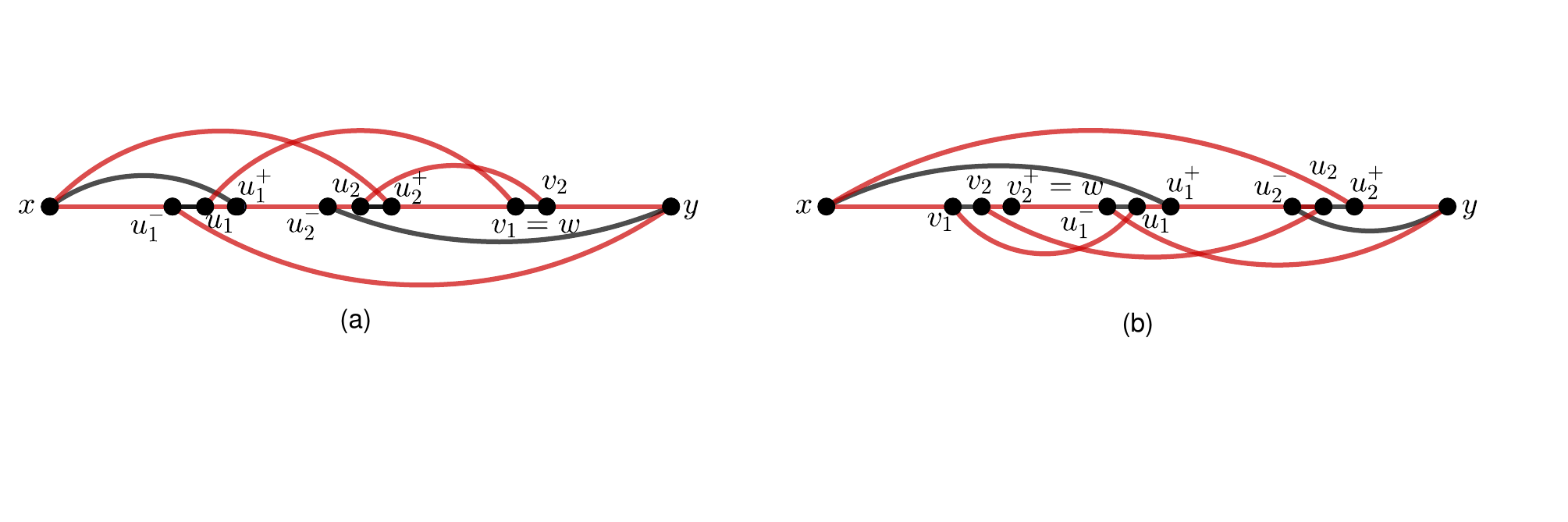}
			\vspace{-2cm}
		\caption{Hamilton cycles in $G$   under the assumption that there exists $w\in (Z\setminus Y) \cup R$ 
			such that there are three distinct vertices $v_1, v_2, v_3\in N(U)\setminus N(U,P)$ for which $g(v_1)=g(v_2)=g(v_3)=w$. }
		\label{f4}
	\end{center}	
\end{figure}

Consider first that  $u_1, u_2$ is on $xPw$.   Then  $g(v_1) =v_1^-$ or $g(v_1)=v_1$ 
and $g(v_2) =v_2^-$ or $g(v_2)=v_2$ by Case 1. 
If $v_1$ is on $xPv_2$, 
 then $g(v_1)=g(v_2)$ implies that $g(v_1)=v_1$, $g(v_2)=v_2^-$,  and $v_2^-=v_1$. 
  However, $xu_2^+Pv_1u_1Pu_2v_2Pyu_1^-Px$ is a Hamilton cycle, see Figure~\ref{f4}(a). 
If $v_2$ is on $xPv_1$, then $g(v_1)=g(v_2)$ implies that $g(v_1)=v^-_1$, $g(v_2)=v_2$,  and $v_1^-=v_2$. 
However, $xu_2^+Pv_2u_2Pu_1v_1Pyu_1^-Px$ is a Hamilton cycle.

Consider then  that  $u_1, u_2$ is on $wPy$.   Then $g(v_1) =v_1^+$ or $g(v_1)=v^{++}_1$ 
and $g(v_2) =v_2^+$ or $g(v_2)=v^{++}_2$ by Case 2. 
If $v_1$ is on $xPv_2$, 
then $g(v_1)=g(v_2)$ implies that $g(v_1)=v^{++}_1$, $g(v_2)=v_2^+$,  and $v_2=v_1^+$. 
However, $xu_2^+Pyu_1^-Pv_2u_2Pu_1v_1Px$ is a Hamilton cycle, see Figure~\ref{f4}(b). 
If $v_2$ is on $xPv_1$, then $g(v_1)=g(v_2)$ implies that $g(v_1)=v^+_1$, $g(v_2)=v^{++}_2$,  and $v_2^+=v_1$. 
However, $xu_2^+Pyu_1^-Pv_1u_1Pu_2v_2Px$ is a Hamilton cycle.

Therefore,  by the arguments above, we have 
\begin{eqnarray*}
	|N(U)| & = &  |N(U)\setminus N(U, P)|+|N(U,P)| \le 2|(Z\setminus Y) \cup R|+|U^-|+|U^+| \le 2t+2|U|. 
\end{eqnarray*}
\qed

Now let $S=N(U)$. 
Then we have $|S|  \le 2t+2|U| $ by Claim~\ref{2.1.7} and  $c(G - S) \geq |U| +1$ by Claim~\ref{2.1.4a}. 
Since $|U| \ge t-1$ by Claim~\ref{2.1.4} and $t\ge 4$, we have  $2|U| <4(|U|-1) \le t(|U|-1)$. 
 Thus  $$\frac{|S|}{c(G - S)} \leq  \frac{2t+2|U|}{1+|U|}< \frac{2t+t(|U|-1)}{2+|U|-1} = t,$$ 
a contradiction to $\tau(G) \ge t$. 

This completes the proof of Theorem~\ref{thm:sum-of-successors}. 
\qqed

\section*{Acknowledgments}
The authors are very grateful to the   two anonymous referees  for their careful reading and valuable comments.

\end{document}